\documentclass[12pt]{article}

\usepackage{graphicx}
\usepackage{pst-poly}     
\usepackage{pst-all}      
\usepackage{amsfonts,amssymb,amsmath,latexsym,bm}

\newtheorem{thm}{Theorem}[section]
\newtheorem{Def}[thm]{Definition}

\newtheorem{lem}[thm]{Lemma}

\newenvironment{pf}[1][Proof]{\noindent\textbf{#1.} }{\hfill\rule{1mm}{2mm}}

\makeatletter \@addtoreset{equation}{section} \makeatother

\begin{document}

\title{A note on the optimal rubbling in ladders and prisms
}
\author
{Zheng-Jiang Xia\footnote{Corresponding author:
xzj@mail.ustc.edu.cn. The research of Zheng-Jiang Xia is supported by Key Projects in Natural Science Research of Anhui Provincial Department of Education (No. KJ2018A0438).},\quad Zhen-Mu Hong\footnote{Email address: zmhong@mail.ustc.edu.cn. The research of Zhen-Mu Hong is supported by NSFC (No.11601002).}
\\
{\small     School of Finance, Anhui University of Finance \& Economics, 233030, Bengbu, China}\\
}

\date{}
\maketitle

\begin{quotation}
\noindent\textbf{Abstract}: A pebbling move on a graph $G$ consists of the removal of two pebbles from one vertex and the placement of one pebble on an adjacent vertex. Rubbling is a version of pebbling where an additional move is allowed, which is also called the strict rubbling move. In this new move, one pebble each is removed from $u$ and $v$ adjacent to a vertex $w$, and one pebble is added on $w$. The optimal rubbling number of a graph $G$ is the smallest number $m$, such that one pebble can be moved to every given vertex from some pebble distribution of $m$ pebbles by a sequence of rubbling moves. In this paper, we give short proofs to determine the rubbling number of cycles and the optimal rubbling number of paths, cycles, ladders, prisms and M$\rm \ddot{o}$bius-ladders.\bigskip

\bigskip

\bigskip
\end{quotation}

\section{Introduction}

Pebbling in graphs was first introduced by Chung\cite{c89}. It has its origin in number theory, and also can be viewed as a model for the transportation of resources, starting from a pebble distribution on the vertices of a connected graph.

Let $G$ be a simple connected graph, we use $V(G)$ and $E(G)$ to denote the vertex set and edge set of $G$, respectively. A pebble distribution $D$ on $G$ is a function $D: V(G)\rightarrow Z$ ($Z$ is the set of nonnegative integers), where $D(v)$ is the number of pebbles on $v$, $|D|=\sum_{v\in V(D)}D(v)$.

A \emph{pebbling move} consists of the removal of two pebbles from a vertex and the placement of one pebble on an adjacent vertex. Let $D$ and $D'$ be two pebble distribution of $G$, we say that $D$ \emph{contains} $D'$ if $D(v)\geq D'(v)$ for all $v\in V(G)$, we say that $D'$ is \emph{reachable} from $D$ if there is some sequence (probably empty) of pebbling moves start from $D$ and resulting in a distribution that contains $D'$.
For a graph $G$, and a vertex $v$, we call $v$ a \emph{root} if the goal is to place pebbles on $v$; If $t$ pebbles can be moved to $v$ from $D$ by a sequence of pebbling moves, then we say that $D$ is $t$-fold \emph{$v$-solvable}, and $v$ is $t$-\emph{reachable} from $D$. If $D$ is $t$-fold $v$-solvable for every vertex $v$, we say that $D$ is $t$-\emph{solvable}.

\emph{The $t$-pebbling number of a graph $G$}, denoted by $f_t(G)$, is the smallest number $m$, such that $t$ pebbles can be moved to every given vertex by pebbling moves from each distribution of $m$ pebbles. While $t=1$, we use $f(G)$ instead of $f_1(G)$, which is called \emph{the pebbling number} of graph $G$. \emph{The optimal pebbling number} $f_{opt}(G)$ of a graph $G$ is the minimum number $m$ for  which there is a pebble distribution $D$ of size $m$ so that every vertex is reachable by a sequence of pebbling moves from $D$.

Rubbling is a version of pebbling where an additional move is allowed, which is also called the strict rubbling move. In this new move, one pebble each is removed from $u$ and $v$ adjacent to a vertex $w$, and one pebble is added on $w$. \emph{The $t$-rubbling number of a graph $G$}, denoted by $\rho_t(G)$, is the smallest number $m$, such that $t$ pebbles can be moved to every given vertex by rubbling moves from each distribution of $m$ pebbles. Similarly, while $t=1$, we use $\rho(G)$ instead of $\rho_1(G)$, which is called \emph{the rubbling number of graph $G$}.
\emph{The optimal rubbling number} $\rho_{opt}(G)$ of a graph $G$ is the minimum number $m$ for  which there is a pebble distribution $D$ of size $m$ so that every vertex is reachable by a sequence of rubbling moves from $D$.

There are many papers about pebbling on graphs, one can view the survey paper\cite{h13} written by G. Hurlbert, rubbling is a new parameter with few results. The basic theory about rubbling and optimal rubbling is developed by~\cite{bs09}, they determined the rubbling number of trees and cycles, the optimal rubbling number of paths and cycles and so on. The rubbling number of complete $m$-ary trees are studied in~\cite{d10}, the rubbling number of caterpillars are given in~\cite{p10}, the optimal rubbling number of ladders, prisms and M$\rm \ddot{o}$bius-ladders are determined in~\cite{kp16}, and in~\cite{ks11}, they give the bounds for the rubbling number of diameter $2$ graphs.

In this paper, we give some new proofs to determine the rubbling number of cycles, the optimal rubbling number of paths, cycles, ladders, prisms and M$\rm \ddot{o}$bius-ladders.



\section{Main Result}
Let $v,w\sim u$, a pebbling move from $v$ to $u$ is the removal of two pebbles from $v$ and addition of one pebble on $u$, denoted by $(v,v\rightarrow u)$. A rubbling move from $\{v,w\}$ to $u$ is the removal of one pebble from $v$ and one pebble from $w$, and addition of one pebble on $u$, denoted by $(v,w\rightarrow u)$.

The following lemma holds clearly.
\begin{lem}{\rm{\cite{bs09}}}\label{lem0}
$2^{d(G)}\leq\rho(G)\leq f(G)$.
\end{lem}

\begin{Def}\rm{\cite{bs09}}\label{def1}
Given a multiset $S$ of rubbling moves on $G$, the \emph{transition digraph} $T(G,S)$ is a directed multigraph whose vertex set is $V(G)$, and each move $(v,w\rightarrow u)$ in $S$ is represented by two directed edges $(v,u)$ and $(w,u)$. The transition digraph of a rubbling sequence $s=(s_1,\ldots,s_n)$ is $T(G,s)=T(G,S)$, where $S=\{s_1,\ldots,s_n\}$ is the multiset of moves in $s$.
\end{Def}

\begin{lem}{\rm(\cite{bs09,ccfh05}, No-Cycle Lemma)}\label{ncl}
Let $S$ be a sequence of rubbling moves on $G$, reaching a distribution $D$. Then there exists a sequence $S^*$ of pebbling moves, reaching a distribution $D^*$, such that

1. on each vertex $v$, $D^*(v)\geq D(v)$;

2. $T(G,S^*)$ does not contain any directed cycles.
\end{lem}

A \emph{thread} in a graph is a path containing vertices of degree 2. By Lemma~\ref{ncl}, we can get

\begin{lem}{\rm\cite{bs09}}\label{lemd}
Let $P$ be a thread in $G$, if vertex $x\notin V(P)$ is reachable from the pebble distribution $D$ by a sequence of rubbling moves, then $x$ is reachable from a rubbling sequence in which there is no strict rubbling move of the form $(v,w\rightarrow u)$ where $u\in V(P)$.
\end{lem}

\begin{Def}\rm\cite{bcc08}
For a given distribution $D$ on $V(G)$, assume the degree of $v$, denoted by $d(v)$, is $2$ and $D(v)\geq 3$. \emph{A smoothing move} from $v$ changes $D$ by remove two pebbles on $v$, and add one pebble on each neighbour of $v$.
\end{Def}

\begin{lem}{\rm\cite{bcc08}\label{lemsmooth}}
Let $D$ be a distribution on a graph $G$ with distinct vertices $u$ and $v$, where $d(v)=2$, $D(v) \geq 3$, and $u$ is $t$-reachable under $D$, then $u$ is $t$-reachable under the
distribution $D'$ obtained by making a smoothing move from $v$.
\end{lem}

A distribution is \emph{smooth}, if it has at most two pebbles on every vertex with degree $2$, a vertex $v$ is \emph{unoccupied} under $D$, if $D(v)=0$.

\begin{lem}{\rm(Smoothing Lemma, \cite{bcc08})\label{smooth}}
If $G$ is a connected $n$-vertex graph, with $n \geq 3$, then $G$
has a smooth optimal distribution with all leaves unoccupied.
\end{lem}

\begin{lem}{\rm\cite{bcc08}}\label{smoothpath}
In a path with a smooth distribution $D$ having at most two pebbles on each endpoint, let $v$ be an unoccupied vertex. If $v$ is an endpoint, then
$v$ is not $2$-reachable under $D$. If $v$ is an inner vertex, then no pebble can be moved out from $v$ without using an edge in both directions.
\end{lem}

\subsection{Rubbling and optimal rubbling in paths and cycles}

\begin{lem}{\rm{\cite{h03,psv95}}}\label{thm2}
The $t$-pebbling numbers of the cycles $C_{2n+1}$ and $C_{2n}$ are
$$f_t(C_{2n+1})=2\left\lfloor\frac{2^{n+1}}{3}\right\rfloor+(t-1)2^n+1,f_t(C_{2n})=2^nt.$$
\end{lem}

In this section, we give a short proof about the rubbling number of cycles and the optimal rubbling number of  paths and cycles, which is determined in~\cite{bs09}.

\begin{thm}{\rm{\cite{bs09}}}\label{thm6}
The rubbling numbers of cycles are
$$\rho(C_{2n})=2^n,\rho(C_{2n+1})=\left\lfloor\frac{7\cdot2^{n-1}-2}{3}\right\rfloor+1.$$
\end{thm}

\begin{pf} For the even cycle $C_{2n}$, by Lemma~\ref{lem0}, we have $2^n\leq \rho(C_{2n})\leq f(C_{2n})=2^n$, over.

For the odd cycle $C_{2n+1}$, assume the target vertex is $v_0$ which is adjacent to $v_1$ and $v_{-1}$. Let $D$ be a pebble distribution on $C_{2n+1}$, and $D(v_0)=0$.
We collapse the vertices $\{v_{-1},v_0,v_1\}$ into one vertex $x$, to get a new graph $C_{2n-1}$ and an induced distribution $D^*$, with $D^*(x)=D(v_{-1})+D(v_0)+D(v_1)$, and $D^*(y)=D(y)$ for $y\neq x$.

By Lemma~\ref{lemd}, we may assume that there is no strict rubbling move of the form $(v,w\rightarrow u)$ where $u\neq v_0$. So
one pebble can be moved to $v_0$ from $D$ on $C_{2n+1}$ by a sequence of rubbling moves $\Leftrightarrow$ two pebbles can be moved to $x$ from $D^*$ on $C_{2n-1}$ by a sequence of pebbling moves. Thus $\rho(C_{2n+1})=f_2(C_{2n-1})=2\left\lfloor\frac{2^n}{3}\right\rfloor+2^{n-1}+1$ (the last equality holds from Lemma~\ref{thm2}). A simple calculation can show that $2\left\lfloor\frac{2^n}{3}\right\rfloor+2^{n-1}+1=\left\lfloor\frac{7\cdot2^{n-1}-2}{3}\right\rfloor+1$, and we are done.
\end{pf}

\begin{lem}{\rm\cite{bcc08}}\label{lemop}
The optimal pebbling numbers of paths and cycles are $$f_{opt}(P_n)=f_{opt}(C_n)=\left\lceil\frac{2n}{3}\right\rceil.$$
\end{lem}

Moreover, we can get the following lemma which is useful later.
\begin{lem}\label{lempath}
Let $D$ be a pebble distribution with $2k$ pebbles on $P_{3k}=v_1\cdots v_{3k}$, if $D$ is solvable by pebbling moves, then $D(v_{3i+2})=2$ for $0\leq i\leq k-1$, and $D(v)=0$ otherwise.
\end{lem}

\begin{pf}
We show it by induction, it holds for $k=1$.

If $D(v_1)\geq1$, then we can move $\left\lfloor D(v_1)\right\rfloor$ pebbles to $P_{3k-1}=v_2\cdots v_{3k}$, then $2k-D(v_1)+\lfloor D(v_1)\rfloor\geq f_{opt}(P_{3k-1})=2k$, thus $D(v_1)=0$. Similarly, $D(v_{3k})=0$.

We make smoothing move on $P_{3k}$ from $D$, to obtain a smooth pebble distribution $D^*$.

Since $|D|=2k<n=3k$, there exist many vertices unoccupied under $D^*$, assume $v_i$ is a vertex unoccupied under $D^*$, Let $L_1=v_1\cdots v_{i}$, and $L_2=v_{i+1}\cdots v_{3k}$.
By Lemma~\ref{smoothpath}, we can not move two pebbles to $v_i$ using only one direction.

If $v_i$ is solvable by $(v_{i-1},v_{i-1}\rightarrow v_i)$, then we can get $D^*|_{L_i}$ is $L_i$-solvable, thus $|D^*|\geq f_{opt}(P_{i})+f_{opt}(P_{3k-i})$. Assume $i=3j+r$ for $0\leq r\leq2$.

If $r=0$, then $|D^*|\geq 2j+2(k-j)=2k$. If $r=1$, then $|D^*|\geq 2j+1+2(k-j)=2k+1$, a contradiction to $|D|=2k$. If $r=2$, then $|D^*|\geq 2j+2+2(k-j-1)+1=2k+1$, a contradiction to $|D|=2k$, thus $i=3j$.
By a similar argument, if $v_i$ is solvable by  $(v_{i+1},v_{i+1}\rightarrow v_i)$, then we can get $i=3j+1$.
Thus we can always partition the path $P_{3k}$ into two paths $L_1=v_1\cdots v_{3j}$ and $L_2=v_{3j+1}\cdots v_{3k}$, and $D^*|_{L_i}$ is solvable in $L_i$ for $i=1,2$, respectively.

By induction, we know that $D^*(v_{3j+2})=2$, and $D^*(v)=0$ otherwise.

Note that making a smooth move on the vertex $v$ leaves at least one pebble on $v$, that means if we can make a smooth move under $D$, then at least two adjacent vertices both having pebbles under $D^*$, a contradiction to $D^*(v_{3j+2})=2$, and $D^*(v)=0$ otherwise. So $D=D^*$, which completes the proof.
\end{pf}

In \cite{bs09}, C. Belford \emph{et al.} determined the optimal rubbling number of paths and cycles, here we give a short proof.

\begin{thm}{\rm\cite{bs09}}\label{thmo1}
The optimal rubbling number of path is $\rho_{opt}(P_n)=\left\lceil\frac{n+1}{2}\right\rceil$, the optimal rubbling number of cycle is $\rho_{opt}(C_n)=\left\lceil\frac{n}{2}\right\rceil$ for $n\geq3$.
\end{thm}

\begin{pf}
\textbf{Upper bound.} For $P_n=v_1v_2\cdots v_n$, let $D$ be a distribution so that $D(v_i)=1$ for $i$ is odd or $i=n$, and $D(v_i)=0$ otherwise, clearly $D$ is solvable, and $|D|=\left\lceil\frac{n+1}{2}\right\rceil$. For $C_n=v_1v_2\cdots v_n$, let $D$ be a distribution so that $D(v_i)=1$ for $i$ is odd , and $D(v_i)=0$ otherwise, clearly $D$ is solvable, and $|D|=\left\lceil\frac{n}{2}\right\rceil$, and we are done.

\textbf{Lower bound.} First we show it for the path $P_n$. We use induction on $n$. By a simple calculation, one can show that it holds for $n\leq3$.

If $n>3$, assume it holds for $P_j$ for all $j<n$. From the upper bound and Lemma~\ref{lemop}, we know that $\rho_{opt}(P_n)<f_{opt}(P_n)$. Thus, let $D$ be an optimal pebble distribution with $\rho_{opt}(P_n)$ pebbles on $P_n$, there must exist a vertex $v_i$, which is reachable only by a strict rubbling move $(v_{i-1},v_{i+1}\rightarrow v_i)$. Let $P_1=v_1\cdots v_{i-1}$, $P_2=v_{i+1}\cdots v_n$ be two subpaths of $P_n$, then by Lemma~\ref{lemd}, we do not need the strict rubbling move $(v_{i-1},v_{i+1}\rightarrow v_i)$ to solve the vertices in $P_i$ for $i=1,2$. That means $D|_{P_i}$ is $v$-solvable for all $v\in P_i$ ($i=1,2$), so $|D|_{P_i}|\geq \rho_{opt}(P_i)$ for $i=1,2$. So $|D|=|D|_{P_1}|+|D|_{P_2}|\geq \left\lceil \frac{i}{2}\right\rceil+\left\lceil\frac{n-i+1}{2}\right\rceil\geq\left\lceil\frac{n+1}{2}\right\rceil$, over.

Now we show it for the cycle $C_n$, similarly, by a simple calculation, one can show that it holds for $n=3$.

If $n>3$, from the upper bound and Lemma~\ref{lemop}, we know that $\rho_{opt}(C_n)<f_{opt}(C_n)$. Thus, let $D$ be an optimal pebble distribution with $\rho_{opt}(C_n)$ pebbles on $C_n$, there must exist a vertex $v_i$, which is reachable only by a strict rubbling move $(v_{i-1},v_{i+1}\rightarrow v_i)$. Then by Lemma~\ref{lemd}, we do not need the strict rubbling move $(v_{i-1},v_{i+1}\rightarrow v_i)$ to solve the vertices in $C_n\backslash v_i$, thus $D$ is solvable on $P_{n-1}=C_n\backslash v_i$, so $|D|\geq\rho_{opt}(P_{n-1})=\left\lceil\frac{n}{2}\right\rceil$, which completes the proof.
\end{pf}

\section{Optimal rubbling in ladders and prisms}

%

%

\begin{Def}\rm\cite{bcc08}
A graph $H$ is a \emph{quotient} of a graph $G$ if the vertices of $H$ correspond to the sets in a partition of $V(G)$,
and distinct vertices of $H$ are adjacent if at least one edge of $G$ has endpoints in the sets
corresponding to both vertices of $H$. In other words, each set in the partition of $V(G)$ \emph{collapses} to a single vertex of $H$.
If $H$ is a quotient of $G$ via the surjective
map $\phi: V (G)\rightarrow V (H)$, and $D$ is a distribution on $G$, then the quotient distribution $D^*$ is
the distribution on $H$ defined by $D^*(u)=\sum_{v\in\phi^{-1}(u)} D(v)$.
\end{Def}

\begin{lem}{\rm( Collapsing Lemma, \cite{bcc08})}\label{collapse}
Let $H$ be a quotient of $G$, $D^*$ in $H$ is induced from $D$ in $G$, then $D$ is $v$-solvable by a sequence of pebbling moves $\Rightarrow$ $D^*$ is $\phi(v)$-solvable by a sequence of pebbling moves.
In particular, $f_{opt}(G)\geq f_{opt}(H)$.
\end{lem}

Similarly, we can get the Collapsing lemma on rubbling.

\begin{lem}\label{rco}
Let $H$ be a quotient of $G$, $D^*$ in $H$ is induced from $D$ in $G$, then $D$ is $v$-solvable by a sequence of rubbling moves $\Rightarrow$ $D^*$ is $\phi(v)$-solvable by a sequence of rubbling moves.
In particular, $\rho_{opt}(G)\geq \rho_{opt}(H)$.
\end{lem}

\begin{pf}
The proof is similar to The proof of Lemma~\ref{collapse} in~\cite{bcc08}.
\end{pf}

%

%

Let $G$ and $H$ be simple connected graphs, we define the \emph{Cartesian product} $G\times H$ to be the graph with vertex
set $V(G\times H)$ and edge set the union of $\{(av, bv)| (a, b)\in E(G),  v\in E(H)\}$  and $\{(ux, uy)|u\in V(G),  (x, y)\in E(H)\}$. We call $P_n\times P_2$ a \emph{ladder} and $C_n\times P_2$ a \emph{prism}. It is clear that a prism can be obtained from a ladder by joining the $4$ endvertices by two edges to form two vertex-disjoint $C_n$ subgraphs. If the four endvertices are joined by two new edges in a switched way to get a $C_{2n}$ subgraph, then a \emph{M$\rm \ddot{o}$bius-ladder} $M_n$ is obtained.

The optimal rubbling numbers of ladders, prisms and M$\rm \ddot{o}$bius-ladders are determined in~\cite{kp16}, here we give new proofs of these results.


\begin{figure}[ht]
\begin{center}
\hspace*{30pt}
\begin{pspicture}(-6,0)(6,4)
\cnode(-5,3){3pt}{u1}\cnode(-3,3){3pt}{u2}\cnode(-1,3){3pt}{u3}\cnode(1,3){3pt}{u4}\cnode(3,3){3pt}{u5}\cnode(5,3){3pt}{u6}
\cnode(-5,1){3pt}{v1}\cnode(-3,1){3pt}{v2}\cnode(-1,1){3pt}{v3}\cnode(1,1){3pt}{v4}\cnode(3,1){3pt}{v5}\cnode(5,1){3pt}{v6}
\ncline{u1}{u2}\ncline{u2}{u3}\ncline{u3}{u4}\ncline[linestyle=dashed,dash=3pt 2pt]{u4}{u5}\ncline{u5}{u6}
\ncline{v1}{v2}\ncline{v2}{v3}\ncline{v3}{v4}\ncline[linestyle=dashed,dash=3pt 2pt]{v4}{v5}\ncline{v5}{v6}
\ncline{u1}{v1}\ncline{u2}{v2}\ncline{u3}{v3}\ncline{u4}{v4}\ncline{u5}{v5}\ncline{u6}{v6}
\ncline[linestyle=dashed,dash=4pt 3pt]{u2}{v1}\ncline[linestyle=dashed,dash=4pt 3pt]{u3}{v2}\ncline[linestyle=dashed,dash=4pt 3pt]{u4}{v3}
\ncline[linestyle=dashed,dash=4pt 3pt]{u6}{v5}

\rput(-5,3.35){$v_{_1}x$}\rput(-3,3.35){$v_{_2}x$}\rput(-1,3.35){$v_{_3}x$}\rput(1,3.35){$v_{_4}x$}\rput(3,3.35){$v_{_{n-1}}x$}
\rput(5,3.35){$v_{_n}x$}
\rput(-5,0.65){$v_{_1}y$}\rput(-3,0.65){$v_{_2}y$}\rput(-1,0.65){$v_{_3}y$}\rput(1,0.65){$v_{_4}y$}\rput(3,0.65){$v_{_{n-1}}y$}
\rput(5,0.65){$v_{_n}y$}

\rput(-5.5,3){$R_{_0}$}\rput(-4.25,2.3){$R_{_1}$}\rput(-2.25,2.3){$R_{_2}$}\rput(-0.25,2.3){$R_{_3}$}\rput(3.75,2.3){$R_{_{n-1}}$}
\rput(5.5,1){$R_{_n}$}
\end{pspicture}
\caption{ $P_n\times P_2$. \label{fig1}
}
%
\hspace*{30pt}
\begin{pspicture}(-6,0)(6,4)
\cnode(-5,3){3pt}{u1}\cnode(-3,3){3pt}{u2}\cnode(-1,3){3pt}{u3}\cnode(1,3){3pt}{u4}\cnode(3,3){3pt}{u5}
\cnode(-5,1){3pt}{v1}\cnode(-3,1){3pt}{v2}\cnode(-1,1){3pt}{v3}\cnode(1,1){3pt}{v4}\cnode(3,1){3pt}{v5}\cnode(5,1){3pt}{v6}
\ncline{u1}{u2}\ncline{u2}{u3}\ncline[linestyle=dashed,dash=3pt 2pt]{u3}{u4}\ncline{u4}{u5}
\ncline{v1}{v2}\ncline{v2}{v3}\ncline[linestyle=dashed,dash=3pt 2pt]{v3}{v4}\ncline{v4}{v5}\ncline{v5}{v6}
\ncline{u1}{v1}\ncline{u2}{v2}\ncline{u3}{v3}\ncline{u4}{v4}\ncline{u5}{v5}

\rput(-5,3.35){$v_{_1}x$}\rput(-3,3.35){$v_{_2}x$}\rput(-1,3.35){$v_{_3}x$}\rput(1,3.35){$v_{_{3k+1}}x$}\rput(3,3.35){$v_{_{3k+2}}x$}
\rput(-5,0.65){$v_{_1}y$}\rput(-3,0.65){$v_{_2}y$}\rput(-1,0.65){$v_{_3}y$}\rput(1,0.65){$v_{_{3k+1}}y$}\rput(3,0.65){$v_{_{3k+2}}y$}
\rput(5,0.65){$v_{_{3k+3}}y$}

\rput(-5.5,2){$R_{_1}$}\rput(-3.5,2){$R_{_2}$}\rput(-1.5,2){$R_{_3}$}\rput(0.5,2){$R_{_{3k+1}}$}\rput(2.5,2){$R_{_{3k+2}}$}
\rput(5.5,1.35){$R_{_{3k+3}}$}
\end{pspicture}
\caption{$H_{3k+2}$. \label{fig2}
}
\end{center}
\end{figure}
%
%
%

Let $P_n=v_1v_2\cdots v_n$ and $P_2=xy$. The vertices of $P_n\times P_2$ are denoted by $v_ix$ and $v_iy$ for $1\leq i\leq n$. Let $H_n=(P _{n+1}\times P_2)\backslash v_{n+1}x$. Then we have

\begin{thm}\label{main}
Let $n=3k+r$, $0\leq r\leq2$,
the optimal rubbling number of the ladder $P_n\times P_2$~{\rm\cite{kp16}} is
\begin{align*}
\rho_{opt}(P_{3k+r}\times P_2)=\left\{
\begin{array}{ll}
2k+1, & {\rm if}\quad  r=0,\\
2k+2, &  {\rm if}\quad r=1,\\
2k+2, &  {\rm if}\quad r=2.
\end{array}
\right.=\left\lceil\frac{2(n+1)}{3}\right\rceil.
\end{align*}

The optimal rubbling number of $H_n$ is
$$\rho(H_n)=\rho(P_{n-1}\times P_2)+1.$$
\end{thm}

\begin{pf}
\textbf{Upper bound.}
Let $D$ be a distribution on $P_n\times P_2$ or $H_n\times P_2$.

If $r=0$, then $D(v_ix)=1$ for $i\equiv 1~(mod~3)$; $D(v_iy)=1$ for $i\equiv 2~(mod~3)$ and $i\neq n-1$;
$D(v_ny)=2$ and $D(v)=0$ otherwise.

If $r=1$, then $D(v_ix)=1$ for $i\equiv 1~(mod~3)$ and $i\neq n$; $D(v_iy)=1$ for $i\equiv 2~(mod~3)$;
$D(v_ny)=2$ and $D(v)=0$ otherwise.

If $r=2$, for $P_n\times P_2$, then $D(v_ix)=1$ for $i\equiv 1~(mod~3)$; $D(v_iy)=1$ for $i\equiv 2~(mod~3)$;
$D(v_ny)=2$ and $D(v)=0$ otherwise. For $H_n$, we add one more pebble on $v_{n+1}y$.

One can check that $D$ is solvable for all $P_n\times P_2$ and $H_n$, over.

\textbf{Lower bound.}
Note that the diameter of $P_n\times P_2$ is $n$, let $R_i$ be the set of vertices with distance $i$ from $v_1x$ ($1\leq i\leq n$), $R_0=v_1x$, if we collapse $R_i$ into one vertex, then we can get a path $L_{n+1}=R_0\cdots R_{n}$ (see Figure~\ref{fig1}). By Lemma~\ref{lemop}, the optimal pebbling number of $L_{n+1}$ is $f_{opt}(L_{n+1})=\left\lceil\frac{2(n+1)}{3}\right\rceil$. By Lemma~\ref{rco}, $\rho_{opt}(H_n)\geq\rho_{opt}(P_n\times P_2)$.
So we only need to show that $\rho_{opt}(P_n\times P_2)\geq f_{opt}(L_{n+1})$ and $\rho_{opt}(H_{3k+2})\geq 2k+3$.

We use induction on $n$ for both $H_n$ (only need to show while $n\equiv2~(mod~3)$) and $P_n\times P_2$, it holds for $n\leq2$ clearly.

Assume it holds for $h<n$.

\textbf{Case $1.$} First we consider $P_n\times P_2$, then let $D$ be a solvable distribution on $P_n\times P_2$ with $\rho_{opt}(P_n\times P_2)$ pebbles, then we collapse $R_i$ to get a path $L_{n+1}$ with length $n$ and induced pebble distribution $D^*$ on $L_{n+1}$.
By Lemma~\ref{rco}, $D^*$ is solvable in $L_{n+1}$.  we will show that $\rho_{opt}(P_n\times P_2)\geq f_{opt}(L_{n+1})$.

Assume $\rho_{opt}(P_n\times P_2)<f_{opt}(L_{n+1})$, then since $D^*$ is a solvable distribution on $L_{n+1}$, thus there must exist some $R_i$ which can be reachable only by strict rubbling move $(R_{i-1},R_{i+1} \rightarrow R_i)$, let $Q_1=R_0\cdots R_{i-1}$, $Q_2=R_{i+1}\cdots R_n$, let $W_1$ and $W_2$ be the subgraphs of $P_n\times P_2$ which induced $Q_1$ and $Q_2$, respectively. Now we will show that $|D|_{W_j}|\geq \rho_{opt}(W_j)$ for $j=1,2$. Note that $W_1$ is isomorphic to $H_{i-1}$ and $W_2$ is isomorphic to $H_{n-i-1}$.

\textbf{Notation:} For simplify, we will use $R_i$  to denote the vertex subset of $P_n\times P_2$, and a vertex of $L_{n+1}$, for example a rubbling move $(R_{i-1},R_{i+1} \rightarrow R_i)$ in $P_n\times P_2$ means a rubbling move $(v,u\rightarrow w)$ for some $v\in R_{i-1}, u\in R_{i+1}$ and $w\in R_i$. We will use under $D$ or $D^*$ to distinguish the rubbling move in different graphs.

If $D|_{W_j}$ is $W_j$-solvable, then $|D|_{W_j}|\geq \rho_{opt}(W_j)$, we are done.
Assume $D|_{W_j}$ is not $W_j$-solvable for some $j$, without loss of generality, assume  $D|_{W_2}$ is not $W_2$-solvable.
If $i=n-1$, then $W_2=\{v_ny\}$, the strict rubbling move $(R_{n-2},R_{n} \rightarrow R_{n-1})$ means that $D^*(R_{n})\geq1$, that means $D(v_ny)\geq1$, a contradiction to  $D|_{W_2}$ is not $W_2$-solvable. Thus $i<n-1$.

Since $D|_{W_2}$ is not $W_2$-solvable, that means to solve some vertex of $W_2$, we must use the pebbles on $W_1$, since we can move at most one pebble on $R_{i-1}$ from $D|_{W_1}$, thus we must use the strict rubbling move $(R_{i-1},R_{i+1} \rightarrow R_i)$ under $D$.
Consider the vertex $v_iy\in R_i$, it can be reachable by some $(R_{i-1},R_{i+1} \rightarrow v_iy)$, in which the vertex in $R_{i+1}$ must be $v_{i+1}y$ since $d(v_iy,v_{i+2}x)=3>1$. Thus we can move one pebble to $v_{i+1}y$ under $D|_{W_2}$, and we can move at most one pebble to $R_{i+1}$ under $D|_{W_2}$ (otherwise we do not need strict rubbling move to solve $R_i$ under $D^*$).

\textbf{Subcase $1.1$} $D(v_{i+1}y)=0$, thus we can move two pebbles on $v_{i+2}y$ (for $v_{i+1}y$ is just connected to one vertex $v_{i+2}y$ of $W_2$). So $v_{i+2}x$ is reachable under $D$, too.

Note that to solve some vertex in $W_2$, we must use strict rubbling move $(R_{i-1},R_{i+1} \rightarrow R_i)$, if the rubbing move is
$(R_{i-1},v_{i+1}y \rightarrow v_iy)$, since $v_iy$ is just connect only one vertex $v_{i+1}y$ of $W_2$, so to use the pebble on $v_iy$, we must move it to $v_{i+1}y$, here induce a directed cycle $v_{i+1}yv_iv_{i+1}y$, which is useless by Lemma~\ref{ncl}, thus the rubbling move must be
$(R_{i-1},v_{i+1}y \rightarrow v_{i+1}x)$. To use the pebble on $v_{i+1}x$, we must use the rubbling move $(R_{i+2},v_{i+1}x \rightarrow v_{i+2}x)$ (since by Lemma~\ref{ncl}, the target vertex cannot be $v_{i+1}y$). Note that we can move two pebbles to $v_{i+2}y$ under $D|_{W_2}$ (which is used to move one pebble on $v_{i+1}y$), so we can use rubbling move $(v_{i+2}y,v_{i+2}y\rightarrow v_{i+2}x)$ instead, a contradiction to the irreplaceability of the strict rubbling move $(R_{i-1},R_{i+1} \rightarrow R_i)$ to solve $W_2$.

\textbf{Subcase $1.2$ }$D(v_{i+1}y)=1$, then at most one pebble can be moved to $R_{i+2}$ under $D|_{W_2}$, otherwise, one more pebble can be moved from $R_{i+2}$ to $R_{i+1}$, thus there are two pebbles on $R_{i+1}$, and one pebble can be moved to $R_i$ under $D^*$, a contradiction.

By a similar argument, we can show there exist some vertex in $W_2$ so that we must use strict rubbling move $(R_{i-1},v_{i+1}y \rightarrow v_{i+1}x)$, and the next rubbling move must be $(R_{i+2},v_{i+1}x \rightarrow v_{i+2}x)$. Since at most one pebble can be moved to $R_{i+2}$ under $D|_{W_2}$, thus exactly one pebble can be moved to $R_{i+2}$ under $D|_{W_2}$.

(a) The vertex  of $R_{i+2}$ used in $(R_{i+2},v_{i+1}x \rightarrow v_{i+2}x)$ is $v_{i+3}x$, then if continues, the rubbling move must be
$(v_{i+3}y,v_{i+2}x \rightarrow v_{i+2}y)$, but then we can use $(v_{i+1}y,v_{i+3}x \rightarrow v_{i+2}y)$ instead, thus we must use the strict rubbling move $(R_{i-1},R_{i+1} \rightarrow R_i)$ to solve at most one vertex  $v_{i+2}x$. Let $P_{left}=v_{i+3}\cdots v_n$, then the vertices of $P_{left}\times P_2$ is reachable without the strict rubbling moves and the pebbble on $v_{i+1}y$. So $|D|_{W_2}|-1\geq \rho_{opt}(P_{left}\times P_2)$, by induction, $\rho_{opt}(P_{left}\times P_2)+1=\rho_{opt}(W_2)$, over.

(b) The vertex  of $R_{i+2}$ used in $(R_{i+2},v_{i+1}x \rightarrow v_{i+2}x)$ is $v_{i+2}y$, since at most one pebble can be moved to $R_{i+2}$, then we may assume $D(v_{i+2}y)=1$ (if not, then we can move two pebbles to $v_{i+3}y$, then one pebble can be moved to $v_{i+3}x$, this is just (a)). Thus, we use one pebble each on $R_{i-1}$, $v_{i+1}y$ and $v_{i+2}y$ to move one pebble to $v_{i+2}x$. Then we rearrangement the distribution on $P_n\times P_2$ as follows: remove the pebble on $v_{i+1}y$, and add it on $v_{i+2}y$, one can view that we can still solve $R_{i}$ and $R_{i+1}$. Thus in the new distribution, we do not need the strict rubbling move $(R_{i-1},R_{i+1} \rightarrow R_i)$ to solve $W_2$, so $|D|_{W_2}|\geq \rho_{opt}(W_2)$.

Thus, $|D|=|D|_{W_1}|+|D|_{W_2}|\geq\rho_{opt}(W_1)+\rho_{opt}(W_2)=\rho_{opt}(H_{i-1})+\rho_{opt}(H_{n-i-1})$.

Assume that $n=3k+r$ and $i=3j+s$ for $0\leq r,s\leq2$, by induction, we can get  Table~\ref{tab1}, where the lower bound is given by $\rho_{opt}(H_{i-1})+\rho_{opt}(H_{n-i-1})$.

\begin{table}[!h]
\begin{center}
  \begin{tabular}{|c|c|c|c|c|c|c|c|c|c|}
  \hline

     $(r,s)$              &$(0,0)$   &(0,1)   &(0,2)   &$(1,0)$     &(1,1)   \\ \hline
    $\rho_{opt}(H_{i-1})$  &2j+1     & 2j+1   & 2j+2   &2j+1        &2j+1         \\ \hline
    $\rho_{opt}(H_{n-i-1})$&2(k-j)+1 & 2(k-j) & 2(k-j)-1      &2(k-j)+1     & 2(k-j)+1         \\ \hline
    Lower bound            &2k+2        &2k+1     &2k+1       &2k+2           & 2k+2      \\ \hline
    $f_{opt}(L_{n+1})$     &2k+1        &2k+1      &2k+1       &2k+2           &2k+2      \\ \hline
  \end{tabular}
  \vspace{.5cm}

  \begin{tabular}{|c|c|c|c|c|c|c|c|c|c|}
  \hline

     $(r,s)$                &(1,2)     &$(2,0)$     &(2,1)     &(2,2)   \\ \hline
    $\rho_{opt}(H_{i-1})$   & 2j+2     &2j+1        & 2j+1     & 2j+2     \\ \hline
    $\rho_{opt}(H_{n-i-1})$ & 2(k-j)   &2(k-j)+2    & 2(k-j)+1 & 2(k-j)+1     \\ \hline
    Lower bound             & 2k+2        &2k+3       & 2k+2        &2k+3 \\ \hline
    $f_{opt}(L_{n+1})$     &2k+2        &2k+2            &2k+2         &2k+2          \\ \hline
  \end{tabular}
  \caption{Bounds on optimal rubbling of $P_n\times P_2$.}\label{tab1}
\end{center}
\end{table}

From Table~\ref{tab1}, we can find $|D|\geq f_{opt}(L_{n+1})$ , which is a contradiction to the assumption that $\rho_{opt}(P_n\times P_2)<f_{opt}(L_{n+1})$, and we are done.

\textbf{Case $2.$} Now we consider $H_{3k+2}$. Note that the diameter of $H_{3k+2}$ is $3k+3$. let $R_i$ be the vertices set $\{v_ix,v_iy\}$ for $1\leq i\leq 3k+2$, and $R_{3k+3}=\{v_{3k+3}y\}$. If we collapse $R_i$ into one vertex, then we can get a path $L_{3k+3}=R_1\cdots R_{3k+3}$ (see Figure~\ref{fig2}). By Lemma~\ref{lemop}, the optimal pebbling number of $L_{3k+3}$ is $f_{opt}(L_{3k+3})=\left\lceil\frac{2(3k+3)}{3}\right\rceil=2k+2$.
So we only need to show that $\rho_{opt}(H_{3k+2})> f_{opt}(L_{3k+3})$.

Let $D$ be a solvable distribution with $\rho_{opt}(H_{3k+2})$ pebbles on $H_{3k+2}$, $D^*$ is an induced pebble distribution on $L_{3k+3}$,
so $D^*$ is solvable on $L_{3k+3}$ by Lemma~\ref{rco}.

\textbf{Subcase $2.1$} If $\rho_{opt}(H_{3k+2})<f_{opt}(L_{3k+3})$, then since $D^*$ is a solvable distribution on $L_{3k+4}$, thus there must exist some $R_i$ which can be reachable only by strict rubbling move $(R_{i-1},R_{i+1} \rightarrow R_i)$, let $Q_1=R_1\cdots R_{i-1}$, $Q_2=R_{i+1}\cdots R_{3k+3}$, let $W_1$ and $W_2$ be two subgraphs of $H_{3k+2}$ induced $Q_1$ and $Q_2$, respectively.

We will show that $D|_{W_i}$ is $W_i$-solvable. Otherwise, we must use the strict rubbling move $(R_{i-1},R_{i+1} \rightarrow R_i)$ under $D$ to solve one vertex of $W_i$, this move must be one of $(v_{i-1}x,v_{i+1}x\rightarrow v_ix)$ and $(v_{i-1}y,v_{i+1}y\rightarrow v_iy)$. Assume the rubbling move is $(v_{i-1}x,v_{i+1}x\rightarrow v_ix)$, then by Lemma~\ref{ncl}, the target vertex of the next rubbling move using the pebble on $v_ix$ cannot be $v_{i+1}x$ or $v_{i-1}x$, if the target vertex is $v_iy$, we need one more pebble on $R_{i-1}$ or $R_{i+1}$, a contradiction to
the condition that $R_i$ can be reachable only by strict rubbling move $(R_{i-1},R_{i+1} \rightarrow R_i)$.

So $D|_{W_i}$ is $W_i$-solvable, thus $|D|_{W_i}|\geq \rho_{opt}(W_i)$ for $i=1,2$. Note that $W_1$ is isomorphic to $P_{i-1}\times P_2$ and $W_2$ is isomorphic to $H_{3k+2-i}$. Thus, $|D|\geq\rho_{opt}(P_{i-1}\times P_2)+\rho_{opt}(H_{3k+2-i})$.
Assume that $i=3j+s$ for $0\leq s\leq2$, by induction, we can get Table~\ref{tab2}, where the lower bound is given by $\rho_{opt}(P_{i-1}\times P_2)+\rho_{opt}(H_{3k+2-i})$.

\begin{table}[!h]
\begin{center}
  \begin{tabular}{|c|c|c|c|c|c|c|c|c|c|}
  \hline

     $s$                             &$0$     &$1$   &$2$      \\ \hline
    $\rho_{opt}(P_{i-1}\times P_2)$  &2j+1          & 2j+2   &2j+2                \\ \hline
    $\rho_{opt}(H_{3k+2-i})$          &2(k-j)+2    & 2(k-j)+1 & 2(k-j)+1             \\ \hline
    Lower bound                      &2k+3        &2k+3     &2k+3         \\ \hline
    $f_{opt}(L_{3k+3})$               &2k+2        &2k+2      &2k+2         \\ \hline
  \end{tabular}
  \caption{Bounds on optimal rubbling of $H_{3k+2}$.}\label{tab2}
\end{center}
\end{table}

From Table~\ref{tab2}, we can find $\rho_{opt}(H_{3k+2})>f_{opt}(L_{3k+3})$, a contradiction to the assumption $\rho_{opt}(H_{3k+2})<f_{opt}(L_{3k+3})$, over.

\textbf{Subcase $2.2$} If $\rho_{opt}(H_{3k+2})=f_{opt}(L_{3k+3})=2k+2$, then $D^*$ is solvable on $L_{3k+3}$. If there exist a vertex $R_i$ which is reachable only by a strict rubbling move under $D^*$, then by a similar argument of Case $2.1$, we are done.

Thus, we only need to consider the case that all vertices of $L_{3k+3}$ are reachable by only pebbling moves, that is to say a strict rubbling move is not allowed. By Lemma~\ref{lempath}, we have $D^*(R_{3i+2})=2$ for $0\leq i\leq k$, and $D^*(v)=0$ otherwise. Thus, at least one of $\{v_1x,v_1y\}$ is not solvable under $D$, which completes the proof.
\end{pf}

%
%
%
%
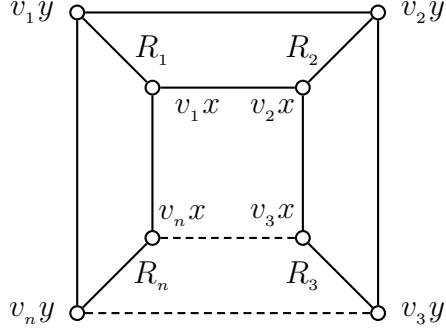
\begin{figure}[ht]
\begin{center}
\hspace*{30pt}
\begin{pspicture}(1,1)(5,5)
\cnode(1,1){3pt}{u1}\cnode(5,1){3pt}{u2}\cnode(5,5){3pt}{u3}\cnode(1,5){3pt}{u4}
\cnode(2,2){3pt}{v1}\cnode(4,2){3pt}{v2}\cnode(4,4){3pt}{v3}\cnode(2,4){3pt}{v4}
\ncline[linestyle=dashed,dash=3pt 2pt]{u1}{u2}\ncline{u2}{u3}\ncline{u3}{u4}\ncline{u4}{u1}
\ncline[linestyle=dashed,dash=3pt 2pt]{v1}{v2}\ncline{v2}{v3}\ncline{v3}{v4}\ncline{v4}{v1}
\ncline{u1}{v1}\ncline{u2}{v2}\ncline{u3}{v3}\ncline{u4}{v4}\ncline{u5}{v5}\ncline{u6}{v6}

\rput(0.4,1){$v_{_n}y$}\rput(5.6,1){$v_{_3}y$}\rput(5.6,5){$v_{_2}y$}\rput(0.4,5){$v_{_1}y$}

\rput(2.4,2.3){$v_{_n}x$}\rput(3.6,2.3){$v_{_3}x$}\rput(3.6,3.7){$v_{_2}x$}\rput(2.6,3.7){$v_{_1}x$}
%
\rput(2,4.5){$R_{_1}$}\rput(4,4.5){$R_{_2}$}\rput(4,1.5){$R_{_3}$}\rput(2,1.5){$R_{_n}$}

\end{pspicture}
\caption{\small ~ $C_n\times P_2$. \label{fig3}
}
\end{center}
\end{figure}

\begin{thm}{\rm\cite{kp16}}
The optimal rubbling number of prism $C_n\times P_2$ is
\begin{align*}
\rho_{opt}(C_{3k+r}\times P_2)=\left\{
\begin{array}{ll}
2k, & {\rm if}\quad  r=0,\\
2k+1, &  {\rm if}\quad r=1,\\
2k+2, &  {\rm if}\quad r=2.
\end{array}
\right.=\left\lceil\frac{2n}{3}\right\rceil.
\end{align*}
except $\rho_{opt}(C_3\times P_2)=3$.
\end{thm}

\begin{pf}
It holds for $n=3$, we only need to show $\rho_{opt}(C_n\times P_2)=\rho_{opt}(P_{n-1}\times P_2)$ for $n\geq4$.

\textbf{Upper bound.} $P_{n-1}\times P_2$ can be viewed as a subgraph of $C_n\times P_2$, it is easy to view that a solvable distribution with $\rho_{opt}(P_{n-1}\times P_2)$ pebbles given in the proof of Theorem~\ref{main} is solvable on $C_n\times P_2$ for $n\geq4$, we are done.

\textbf{Lower bound.} Let $C_n=v_1\cdots v_n$, $P_2=xy$, assume $D$ is a solvable distribution with  $\rho_{opt}(C_n\times P_2)$ pebbles, we collapse each set $\{v_ix,v_iy\}$ into one vertex $R_i$ (see Figure~\ref{fig3}), then we get a cycle $C_n=R_1\cdots R_n$ and an induced distribution $D^*$. By Lemma~\ref{rco}, $D^*$
is solvable on $C_n$. By Lemma~\ref{lemop}, $f_{opt}(C_n)=\left\lceil\frac{2n}{3}\right\rceil$, we only need to show $\rho_{opt}(C_n\times P_2)=f_{opt}(C_n)$.

Assume $\rho_{opt}(C_n\times P_2)<f_{opt}(C_n)$, then there must exist a vertex $R_i$ which is reachable under $D^*$ only by strict rubbling move $(R_{i-1},R_{i+1}\rightarrow R_i)$. By a similar argument of Case 2.1 in the proof of Theorem~\ref{main}, if we remove $\{v_ix,v_iy\}$ from $C_n\times P_2$ (which is just isomorphic to $P_{n-1}\times P_2$), then $D$ is still solvable. So $|D|\geq \rho_{opt}(P_{n-1}\times P_2)=\lceil\frac{2n}{3}\rceil$, which is a contradiction to the assumption $\rho_{opt}(C_n\times P_2)<f_{opt}(C_n)$, and we are done.
\end{pf}

\begin{thm}{\rm\cite{kp16}}
The optimal rubbling number of the M$\rm \ddot{o}$bius ladder $M_n$ is
$$\rho_{opt}(M_n)=
\rho_{opt}(C_{n}\times P_2).$$
\end{thm}

\begin{pf}
The proof is similar to the proof of prisms.
\end{pf}

\end{document}